\documentclass[reqno]{amsart}
\usepackage{hyperref}
\usepackage{amsmath}
\usepackage{amssymb}
\usepackage{amscd}
\usepackage{graphics}
\usepackage{latexsym}
\usepackage{stmaryrd}

\theoremstyle{plain}
\newtheorem{theorem}{Theorem}[section]
\newtheorem{thm}[theorem]{Theorem}
\newtheorem{cor}[theorem]{Corollary}
\newtheorem{lem}[theorem]{Lemma}
\newtheorem{prop}[theorem]{Proposition}
\newtheorem{claim}[theorem]{Claim}

\theoremstyle{definition}
\newtheorem{defn}[theorem]{Definition}

\newtheorem{rmk}[theorem]{Remark}

\newtheorem{notat}[theorem]{Notation}

\theoremstyle{remark}

\newcommand{\ZZ}{\mathbb{Z}}

\newcommand{\AAA}{\mathbb{A}}
\newcommand{\PP}{\mathbb{P}}

\newcommand{\mc}{\mathcal}
\newcommand{\mf}{\mathfrak}

\newcommand{\OO}{\mc{O}}

\newcommand{\SP}{\text{Spec }}

\newcommand{\pflpf}{\text{pr,f,lpf}}

\newcommand{\lt}{\left}
\newcommand{\rt}{\right}

\newsavebox{\sembox}
\newlength{\semwidth}
\newlength{\boxwidth}

\newcommand{\Sem}[1]{%
\sbox{\sembox}{\ensuremath{#1}}%
\settowidth{\semwidth}{\usebox{\sembox}}%
\sbox{\sembox}{\ensuremath{\left[\usebox{\sembox}\right]}}%
\settowidth{\boxwidth}{\usebox{\sembox}}%
\addtolength{\boxwidth}{-\semwidth}%
\left[\hspace{-0.3\boxwidth}%
\usebox{\sembox}%
\hspace{-0.3\boxwidth}\right]%
}

\newsavebox{\semrbox}
\newlength{\semrwidth}
\newlength{\boxrwidth}

\title[Moduli spaces]{Artin's axioms, composition and moduli spaces} 
\author[Starr]{Jason Starr} 
\begin{document}


\begin{abstract}
We prove Artin's axioms for algebraicity of a stack are compatible
with composition of 1-morphisms.  Consequently, some natural stacks
are algebraic.  One of these is a common generalization of Vistoli's
Hilbert stack and the stack of branchvarieties defined by Alexeev and
Knutson.
\end{abstract}


\maketitle


\section{Introduction} \label{sec-intro}

Many moduli functors in algebraic geometry, properly interpreted, are
\emph{algebraic stacks}, also called \emph{Artin stacks} (please note,
following Artin, we do not assume diagonal morphisms are
quasi-compact).  In ~\cite{Versal}, Artin gave axioms for algebraicity
involving deformation-obstruction theory and compatibility with
completion.

There exist natural stacks $\mc{Y}$ where the completion axiom fails
but all other axioms hold.  Sometimes there exists a stack $\mc{X}$
for which the completion axiom holds and a 1-morphism
$f:\mc{X}\rightarrow \mc{Y}$ representable by algebraic stacks.
Intuitively $\mc{X}$ should satisfy all Artin's axioms, and thus be
algebraic.  In other words, Artin's axioms should be compatible with
\emph{composition} of 1-morphisms of stacks.  The difficulty is that,
given a relative obstruction theory for $f$ and an obstruction theory
for $\mc{Y}$, there may exist no ``extension'' obstruction theory for
$\mc{X}$.

Existence of an extension obstruction theory is circumvented using
Propositions ~\ref{prop-comp1} and ~\ref{prop-comp2}.  The main result
is the following version of ~\cite[Theorem 5.3]{Versal}.

\begin{prop} \label{prop-main}
Let $S$ be an excellent scheme and let $f:\mc{X} \rightarrow \mc{Y}$
be a 1-morphism of limit preserving stacks in groupoids over
$(\text{Aff}/S)$ for the \'etale topology.  Let $O_\mc{Y}$ be an
obstruction theory for $\mc{Y}$ and let $O_f$ be a relative
obstruction theory for $f$.  The stack $\mc{X}$ is algebraic if,
\begin{enumerate}
\item[(1)]
Conditions ~\cite[(S1,2)]{Versal} hold for deformations and
automorphisms of $\mc{Y}$ and $\mc{X}$ (or equivalently, $\mc{Y}$ and
$f$).
\item[(2)]
For any complete local $\OO_S$-algebra $\widehat{A}$ with residue
field of finite type over $S$, the canonical map
$$
\mc{X}(\SP \widehat{A}) \rightarrow \varprojlim \mc{X}(\SP
\widehat{A}/\mf{m}^n)
$$
is faithful, and has a dense image, i.e., the projection to
$\mc{X}(\SP \widehat{A}/\mf{m}^n)$ is essentially surjective for every $n$.
\item[(3)]
Automorphisms, deformations and obstructions of $\mc{Y}$ and $f$
satisfy the conditions in Notations ~\ref{notat-41a} and
~\ref{notat-41}.
\item[(4)]
If the object $a_0$ of $\mc{X}(\SP A_0)$ is algebraic, and if $\phi$
is an automorphism of $a_0$ inducing the identity in $\mc{X}(\SP
k(y))$ for a dense set of finite type points $y$ of $\SP A_0$, then
$\phi$ equals $\text{Id}_{a_0}$ on a non-empty open subset of $\SP
A_0$.
\end{enumerate}
\end{prop}

A consequence is algebraicity of some natural stacks.

\begin{prop} \label{prop-main2}
Let $S$ be an excellent scheme and let $\mc{Y}$ be a limit preserving
algebraic stack over $(\text{Aff}/S)$ with finite diagonal.  The stack
$\mc{H}$ parametrizing triples $(X,L,g)$ of a proper algebraic space
$X$, a 1-morphism $g:X\rightarrow \mc{H}$, and an invertible,
$g$-ample $\OO_X$-module $L$ is a limit preserving algebraic stack
over $(\text{Aff}/S)$ with quasi-compact, separated diagonal.
\end{prop}

This stack $\mc{H}$ is a common generalization of Vistoli's
\emph{Hilbert stack}, ~\cite{VHilb}, and the stack of
\emph{branchvarieties} defined by Alexeev and Knutson, ~\cite{AlKn}.
The proof of Proposition ~\ref{prop-main2} gives a new proof of
algebraicity in each of these special cases.

\medskip\noindent
\textbf{Acknowledgments.}  I thank Martin Olsson for pointing out
~\cite[Appendix A]{OK3}.  I thank Valery Alexeev and Allen Knutson for
pointing out that Proposition ~\ref{prop-main2} is not a formal
consequence of Vistoli's Hilbert stack.


\section{Artin's axioms and obstruction theory} \label{sec-aa}

All hypotheses regarding obstructions in ~\cite{Versal} trace back to
the proofs of ~\cite[Proposition 4.3 and Theorem 4.4]{Versal}.  This
section describes lifting properties, why the lifting properties are
compatible with composition, how Artin uses these lifting properties
to prove openness of versality, and finally, how Artin uses
obstruction theory to prove the lifting properties.

\subsection{Relative deformation situations} \label{ssec-not}

Relative deformation situations and obstruction theories are studied
in ~\cite[Appendix A]{OK3}.  Following are the basic definitions.

Let $S$ be a locally Noetherian algebraic space (for later
applications, it will be necessary that $S$ is excellent).  An
\emph{infinitesimal extension} is a surjective homomorphism of
Noetherian $\OO_S$-algebras with nilpotent kernel, $A'\xrightarrow{q}
A$.  An \emph{extension pair} is a pair of infinitesimal extensions
$(A'\xrightarrow{q} A \xrightarrow{q_0} A_0)$ such that $A_0$ is
reduced and the kernel $M$ of $q$ is annihilated by the kernel of
$q_0\circ q$.

A \emph{morphism of infinitesimal extensions} $(u',u)$ is a Cartesian
diagram of $\OO_S$-algebras,
$$
\begin{CD}
A' @> q_A >> A \\
@V u' VV  @VV u V \\
B' @> q_B >> B
\end{CD}
$$
whose rows are infinitesimal extensions.  \emph{Morphisms of extension
pairs} $(u',u,u_0)$ are defined analogously.

Let $\mc{X}$ be a stacks in groupoids over $(\text{Aff}/S)$ for the
\'etale topology.  An \emph{infinitesimal extension over $\mc{X}$} is
a datum $(q,a)$ of a reduced infinitesimal extension
$A\xrightarrow{q_0} A_0$ and an object $a$ of $\mc{X}(\SP A)$.  A
\emph{morphism of infinitesimal extensions over $\mc{X}$},
$$
(u,u_0,\phi):(A\xrightarrow{q_{A,0}} A_0,a) \rightarrow
(B\xrightarrow{q_{B,0}} B_0,b),
$$
is a morphism $(u,u_0)$ of infinitesimal extensions together with a
morphism $\phi:b \rightarrow a$ in $\mc{X}$ mapping to $u^*:\SP
B\rightarrow \SP A$.

\begin{defn} \label{defn-ds}
A \emph{deformation situation over $\mc{X}$} is a datum
$(A'\xrightarrow{q} A \xrightarrow{q_0}A_0,a)$ of an extension pair
$(q,q_0)$ and an object $a$ of $\mc{X}(\SP A)$.

A \emph{morphism of deformation situations over $\mc{X}$} is a datum
$(u',u,u_0,\phi)$ of a morphism $(u',u,u_0)$ of extension pairs
together with a morphism $\phi:b\rightarrow a$ in $\mc{X}$ mapping to
$u^*$.
\end{defn}
 
By the axioms for a stack, given a deformation situation $(q,q_0,a)$
and a morphism of extension pairs $(u',u,u_0)$, there exists a
morphism $\phi$ so that $(u',u,u_0,\phi)$ is a morphism of deformation
situations over $\mc{X}$.  A \emph{clivage normalis\'e} determines a
choice of $\phi$, cf. ~\cite[D\'efinition VI.7.1]{SGA1}.  From this
point on, a clivage normalis\'e is assumed given.  The image of
$(u',u,u_0,\phi)$ is called the \emph{base change} of $(q,q_0,a)$ by
$(u',u,u_0)$.

Let $f:\mc{X}\rightarrow \mc{Y}$ be a 1-morphism of stacks in
groupoids over $(\text{Aff}/S)$.  An \emph{infinitesimal extension
over $f$} is a datum $(\widetilde{q},a)$ of a morphism $\widetilde{q}$
of $\mc{Y}$ mapping to a reduced infinitesimal extension
$A\xrightarrow{q_0} A_0$ and an object $a$ of $\mc{X}(\SP A)$ mapping
to the target of $\widetilde{q}$.  Equivalently, it is a 1-morphism
$\SP A \rightarrow \mc{Y}$ and an infinitesimal extension over the
2-fibered product $\SP A \times_{\mc{Y}} \mc{X}$.  Morphisms of
infinitesimal extensions over $f$, deformation situations over $f$,
morphisms of deformation situations over $f$, and base change are
defined analogously.

If $\mc{Y}$ satisfies the Schlessinger-Rim criterion
~\cite[(S1)]{Versal}, there are well-defined relative analogues of the
Schlessinger-Rim criterion for $f$.

\begin{lem} ~\cite[\S A.14]{OK3} \label{lem-sorites}
Let $\mc{Y}$ be a stack in groupoids over $(\text{Aff}/S)$ satisfying
~\cite[(S1)]{Versal}.  Let $f:\mc{X}\rightarrow \mc{Y}$ be a
1-morphism of stacks in groupoids over $(\text{Aff}/S)$.
\begin{enumerate}
\item[(i)] 
The stack in groupoids $\mc{X}$ satisfies ~\cite[(S1)]{Versal} if and
only if $f$ satisfies the relative analogue of ~\cite[(S1)]{Versal}.
\item[(ii)]
Let $g:\mc{Z}\rightarrow \mc{Y}$ be a 1-morphism of stacks in
groupoids.  If $\mc{X}$ and $\mc{Z}$ satisfy ~\cite[(S1)]{Versal},
then also the 2-fibered product $\mc{X}\times_{\mc{Y}}\mc{Z}$
satisfies ~\cite[(S1)]{Versal}.
\end{enumerate}
\end{lem}

\begin{proof}
This is largely verified in ~\cite[\S A.14]{OK3}.  The details are
left to the reader.
\end{proof}

Let $f:\mc{X}\rightarrow \mc{Y}$ be a 1-morphism of stacks in
groupoids over $(\text{Aff}/S)$, each of which satisfies
~\cite[(S1)]{Versal}.  For each reduced $\OO_S$-algebra $A_0$ and
object $a_0$ of $\mc{X}(\SP A_0)$, ~\cite[\S A.15]{OK3} gives a
natural 7-term exact sequence of automorphism and deformation groups.
This implies the following.

\begin{lem} ~\cite[\S A.15]{OK3} \label{lem-sorites2}
\begin{enumerate}
\item[(i)] 
Assuming $\mc{Y}$ satisfies ~\cite[(S2)]{Versal}, $\mc{X}$ satisfies
~\cite[(S2)]{Versal} if and only if $f$ satisfies the relative
analogue of ~\cite[(S2)]{Versal}.
\item[(ii)]
If $g:\mc{Z}\rightarrow \mc{Y}$ is a 1-morphism of stacks in groupoids
such that $\mc{X},\mc{Y}$ and $\mc{Z}$ each satisfy
~\cite[(S1,2)]{Versal}, then also $\mc{X}\times_{\mc{Y}} \mc{Z}$
satisfies ~\cite[(S1,2)]{Versal}.
\end{enumerate}
\end{lem}

\begin{proof}
This is largely verified in ~\cite[\S A.15]{OK3}.  The details are
left to the reader.
\end{proof}

\begin{notat}~\cite[(4.1)]{Versal}, ~\cite[\S A.11]{OK3} \label{notat-41a}
There are relative analogues of the conditions on automorphisms and
deformations. (The relative analogues of conditions on obstructions
are stated in Notation ~\ref{notat-41}.)
\begin{enumerate}
\item[(4.1.i)]
The functors, $\text{Aut}_f$, resp. $D_f$, are compatible with \'etale
localization: For every morphism
$(\widetilde{u},\widetilde{u}_0,\phi)$ of infinitesimal extensions
over $f$, if $A_0$ is a finite-type $\OO_S$-algebra and $u$ is
\'etale, then the following associated natural transformations of
functors are isomorphisms,
$$
\text{Aut}_{f,a_0}(A_0+M)\otimes_{A_0} B_0 \rightarrow
\text{Aut}_{f,b_0}(B_0+M\otimes_{A_0}B_0), 
$$
$$
D_{f,a_0}(M)\otimes_{A_0} B_0 \rightarrow
D_{f,b_0}(M\otimes_{A_0} B_0). 
$$
\item[(4.1.ii)]
The functors $\text{Aut}_f$ and $D_f$ are compatible with completions:
For every finite-type $\OO_S$-algebra $A_0$ and every maximal ideal
$\mathfrak{m}$ of $A_0$, the following natural maps are isomorphisms,
$$
\text{Aut}_{f,a_0}(A_0+M)\otimes_{A_0} \widehat{A}_0 \rightarrow \varprojlim
\text{Aut}_{f,a_0}((A_0+M)/\mathfrak{m}^n),
$$
$$
D_{f,a_0}(M)\otimes_{A_0}\widehat{A}_0 \rightarrow \varprojlim
D_{f,a_0}(M/\mathfrak{m}^nM).
$$
\item[(4.1.iii)]
For every infinitesimal extension $(\widetilde{q}_0,a)$ over $f$ with
$A$ a finite type $\OO_S$-algebra, there is an open dense set of
points of finite type $p\in \SP A_0$ so that the following maps are
isomorphisms,
$$
\text{Aut}_{f,a_0}(A_0+M)\otimes_{A_0}k(p) \rightarrow
\text{Aut}_{f,a_0}(k(p) + M\otimes_{A_0}k(p)),
$$
$$
D_{f,a_0}(M)\otimes_{A_0} k(p) \rightarrow
D_{f,a_0}(M\otimes_{A_0}k(p)).
$$
\end{enumerate}
\end{notat}

\begin{lem} \label{lem-sorites3}
With notation as in Lemma~\ref{lem-sorites} and
Lemma~\ref{lem-sorites2}, assume $\mc{X}$ and $\mc{Y}$ satisfy
~\cite[(S1,2)]{Versal} and assume automorphisms, respectively
deformations, of $\mc{Y}$ satisfy ~\cite[(4.1)]{Versal}.
Automorphisms, resp. deformations, of $\mc{X}$ satisfy
~\cite[(4.1)]{Versal} if and only if automorphisms,
resp. deformations, of $f$ satisfy ~\ref{notat-41}.  Also, if
$f:\mc{X}\rightarrow \mc{Y}$ and $g:\mc{Z}\rightarrow \mc{Y}$ are
1-morphisms of stacks in groupoids over $(\text{Aff}/S)$ whose
automorphisms, resp. deformations, satisfy
~\cite[(S1,2),(4.1)]{Versal}, then also the automorphisms, resp.
deformations, of $\mc{X}\times_{\mc{Y}}\mc{Z}$ satisfy
~\cite[(S1,2),(4.1)]{Versal}.
\end{lem}

\begin{proof}
This also largely follows from ~\cite[\S A.15]{OK3}.  The details are
left to the reader.
\end{proof}

\subsection{Lifting properties} \label{ssec-lift}
Let $(B'\xrightarrow{q} B\xrightarrow{q_0} B_0)$ be an extension pair
and let $(\widetilde{q},\widetilde{q}_0,b_{\mc{X}})$ be a deformation
situation over $f$ mapping to $(q,q_0)$.  The deformation situation is
\emph{algebraic}, resp. \emph{reduced}, \emph{integral}, if $B_0$ is
finite type, resp. reduced, integral.

\medskip\noindent
\textbf{Localization.} The lifting properties involve localization.
Let $q:B'\rightarrow B$ be an infinitesimal extension.  Images and
inverse images of multiplicative systems under $q$ are again
multiplicative.  Let $S'$ be a multiplicative subset of $B'$ and let
$S'_q$ be the multiplicative subset $q^{-1}(q(S'))$.  The localization
of $B'$ with respect to $S'$ equals the localization of $B'$ with
respect to $S'_q$.  Moreover, $S^{-1}B'\rightarrow q(S)^{-1}B$ is an
infinitesimal extension, and the associated graded pieces of the
kernel are the localizations $q(S)^{-1}(N^i/N^{i+1})$.  Therefore, the
localizations of $B'$ are in 1-to-1 bijection with the localizations
of $B$.  Moreover, a localization of $B$ is finitely generated if and
only if the associated localization of $B'$ is finitely generated.

Let $(B'\xrightarrow{\widetilde{q}_B}
B\xrightarrow{\widetilde{q}_{B,0}} B_0, b_{\mc{X}})$ be an integral
deformation situation.  Let $A_0$ be the fraction field of $B_0$, and
let $A$ and $A'$ be the associated localizations of $B$ and $B'$
respectively.  Let $A' \xrightarrow{q_A} A \xrightarrow{q_{A,0}} A_0$
denote the associated extension pair.  Denote by
$(a_{\mc{Y}}'\xrightarrow{\widetilde{q}_A} a_{\mc{Y}}
\xrightarrow{q_{A,0}} a_{\mc{Y},0},a_{\mc{X}})$ the base change over
$A'$, etc., of the deformation situation
$(\widetilde{q}_B,\widetilde{q}_{B,0},b)$.

\begin{defn} \label{defn-lift}
A \emph{generic lift} of the deformation situation
$(\widetilde{q}_B,\widetilde{q}_{B,0},b_{\mc{X}})$ is a morphism
$a_{\mc{X}}' \xrightarrow{r_A} a_{\mc{X}}$ in $\mc{X}$ over
$\widetilde{q}_A$.  An \emph{integral lift} of the generic lift is a
morphism $b'_{\mc{X},\text{new}} \xrightarrow{r_B} b_{\mc{X}}$ in
$\mc{X}$ satisfying the following conditions.
\begin{enumerate}
\item[(i)]
The image of $r_B$ in $(\text{Aff}/S)$ is an infinitesimal extension
$B'_\text{new} \xrightarrow{q_{B,\text{new}}} B$ sitting between $B'$
and $q_{A}^{-1}(B)$.
\item[(ii)]
The image $\widetilde{q}_{B,\text{new}}$ of $r_B$ in $\mc{Y}$ is the
base-change of $\widetilde{q}_B$.
\item[(iii)]
And the base-change of $r_B$ by $(q_{B,\text{new}},q_{B,0})\rightarrow
(q_A,q_{A,0})$ is $r_A$.
\end{enumerate}
\end{defn} 

\medskip\noindent
\textbf{Pushouts of $M$.}  Given a deformation situation
$(\widetilde{q},\widetilde{q}_{0},b_{\mc{X}})$ over $f$ mapping to an
extension pair $(B'\xrightarrow{q} B\xrightarrow{q_0} B_0)$, for every
surjection of $B_0$-modules, $M\rightarrow N$, there is a surjection
of $\OO_S$-algebras $B'\xrightarrow{u} B'_N$ whose kernel is the
kernel of $M\rightarrow N$.  This gives a morphism of deformation
situations,
$$
(\widetilde{u}',\text{Id},\text{Id},\text{Id}):
(\widetilde{q},\widetilde{q}_{0},b_{\mc{X}}) \rightarrow
(\widetilde{q}_{N}, \widetilde{q}_{0},b_{\mc{X}}),
$$ 
such that the image of $\widetilde{q}_{N}$ is the extension
$B'_N\rightarrow B$ and the image of $\widetilde{u}'$ is $u$.

Let $\mf{a}$ be a radical ideal in $B_0$.  Let $R$ denote the
semilocalization of $B_0$ at the generic points of $\mf{a}$.  Assume
$M$ is a finite type $B_0/\mf{a}$-module.  A localization $A_0$ of
$B_0$ is \emph{$\mf{a}$-generic} if $\SP A_0 \cap \SP(B_0/\mf{a})$ is
dense in $\SP (B_0/\mf{a})$, i.e., $A_0$ is isomorphic to a
$B_0$-subalgebra of $R$.  An \emph{$\mf{a}$-generic quotient} of $M$
is a pair $(A_0,N)$ of an $\mf{a}$-generic localization $A_0$ and a
surjection $M\otimes_{B_0} A_0 \rightarrow N$.  It is \emph{finite
type} if $A_0$ is a finite type $B_0$-algebra.  It is
\emph{projective} if $N$ is a projective $A_0/\mf{a}A_0$-module.  It
is \emph{extending} if for the associated deformation situation
$(\widetilde{q}_{A,N},\widetilde{q}_{A,0},a_{\mc{X}})$ over
$(A'_N,A,A_0)$, there is a lifting of $\widetilde{q}_{A,N}$ to
$\mc{X}$.

\begin{defn} \label{defn-extend}
Let $(\widetilde{q},\widetilde{q}_0,b_{\mc{X}})$ be a deformation
situation over $f$ whose kernel $M$ is a finite type
$B_0/\mathfrak{a}$-module for a radical ideal $\mathfrak{a}$.  A
\emph{generic extender} is a finite type, projective, extending,
$\mf{a}$-generic quotient $(A_0,N)$ such that for every finite type,
projective, $\mf{a}$-generic quotient $(C_0,P)$, the quotient
$P\otimes_{C_0} R$ factors through the quotient $N\otimes_{A_0} R$.

A generic extender is \emph{compatible with \'etale extension} if for
every \'etale homomorphism $v: B'\rightarrow B'_{\text{\'et}}$, the
pair
$(A_0\otimes_{B_0}B_{\text{\'et},0},N\otimes_{B_0}B_{\text{\'et},0})$
is a $v(\mf{a})$-generic extender for the base change of the
deformation situation to $B'_{\text{\'et}}$.

A generic extender is \emph{compatible with closed points} if for
every closed point $A_0\rightarrow k(y)$ of $\SP A_0$ and every
surjection $M\otimes_{B_0} k(y) \rightarrow N_y$, there exists a lift
$\widetilde{q}_{\mc{X},N_y}$ of the base change of $\widetilde{q}$ to
$A'_{N_y}$ if and only if $N_y$ is a quotient of $N\otimes_{A_0}
k(y)$.
\end{defn}

\begin{lem} \label{lem-extend}
Assume $\mc{X}$ and $\mc{Y}$ each satisfy ~\cite[(S1)]{Versal}.  Then
for every deformation situation and radical ideal $\mf{a}$ there
exists a generic extender.
\end{lem}

\begin{proof}
For every finite type, projective, extending, $\mf{a}$-generic
quotient $(C_0,P)$, consider the quotient $M\otimes_{B_0}R\rightarrow
P\otimes_{C_0} R$.  This system of quotients has an inverse limit.
Because $M\otimes_{B_0} R$ has finite length, the inverse limit is
equal to the inverse limit of a finite subsystem.  Thus, it suffices
to prove the following.  For every pair of finite type, projective,
extending, $\mf{a}$-generic quotients $(A_{0,1},N_1)$ and
$(A_{0,2},N_2)$, there exists a finite type, projective, extending,
$\mf{a}$-generic quotient $(A_0,N)$ such that both quotients
$N_1\otimes_{A_{0,1}} R$ and $N_2\otimes_{A_{0,2}} R$ factor through
$N\otimes_{A_0} R$.

Replace $B_0$ by the finite type, $\mf{a}$-generic localization
$A_{0,1}\otimes_{B_0} A_{0,2}$ and replace $B$ and $B'$ by the
associated localizations.  Consider the induced map $M\rightarrow
N_1\oplus N_2$.  After a further finite type, $\mf{a}$-generic
localization, the cokernel $P$ is a projective $B_0/\mf{a}$-module.
Denote by $N$ the image of $M$ in $N_1\oplus N_2$.  Of course $B'_N$
equals the fiber product $B'_{N_1}\times_{B'_P} B'_{N_2}$.

Let $D$ denote the functor $D_{f,\widetilde{q}_0,b_{\mc{X},0}}$.  For
$i=1,2$, let $\widetilde{q}_{\mc{X},i}:b_{\mc{X},i} \rightarrow
b_{\mc{X}}$ be a lift to $\mc{X}$ of the base change of
$\widetilde{q}$ over $B'_{N_i}$.  The base change of $b_{\mc{X},i}$
over $B'_P$ is a lift to $\mc{X}$ of the base change of
$\widetilde{q}$ over $B'_P$.  Therefore the 2 base changes differ by
an element $d$ in $D(P)$.

Because $N_1\oplus N_2\rightarrow P$ is a surjective map of
$B_0/\mf{a}$-modules whose image is a projective, it is split.
Therefore $D(N_1\oplus N_2)\rightarrow D(P)$ is surjective.  So there
exist elements $d_i$ of $D(N_i)$ for $i=1,2$ such that $d$ is the
image of $(d_1,d_2)$.  After translating $b_{\mc{X},1}$,
resp. $b_{\mc{X},2}$, by $d_1$, resp. $d_2$, the base changes to
$B'_P$ agree.  Therefore by ~\cite[(S1)]{Versal}, there exists an
element $b_{\mc{X}}$ over $B'_N$ whose base change to $B'_{N_i}$
equals $b_{\mc{X},i}$ for $i=1,2$.
\end{proof}

\subsection{Lifting properties and compositions} \label{ssect-comp}
Let $e:\mc{W}\rightarrow \mc{X}$ and $f:\mc{X}\rightarrow \mc{Y}$ be
1-morphisms of stacks in groupoids over $(\text{Aff}/S)$.

\begin{prop} \label{prop-comp1}
If for both $e$ and $f$ generic lifts of integral, algebraic
deformation situations have integral lifts, then the same holds for
$f\circ e$.
\end{prop}

\begin{proof}
Given an integral, algebraic deformation situation
$(\widetilde{q}_{\mc{Y},B},\widetilde{q}_{\mc{Y},B,0},b_{\mc{W}})$
over $f\circ e$ and a generic lift $a_{\mc{W}}'\rightarrow
a_{\mc{W}}$, then $e(a_{\mc{W}}')\rightarrow e(a_{\mc{W}})$ is a
generic lift of the integral, algebraic deformation situation
$(\widetilde{q}_{\mc{Y},B}, \widetilde{q}_{\mc{Y},B,0},e(b_{\mc{W}}))$
over $f$.  By hypothesis, this has an integral lift $r_f$.  After
replacing $B'$, this integral lift gives a deformation situation
$(r_f,\widetilde{q}_{\mc{X},B,0},b_{\mc{W}})$ over $e$.  And
$a_{\mc{W}}'\rightarrow a_{\mc{W}}$ is a generic lift.  By hypothesis,
there exists an integral lift $r_e$.  This is also an integral lift
for the original generic lift over $f\circ e$.
\end{proof}

\begin{prop} \label{prop-comp2}
Assume each of $\mc{W}$, $\mc{X}$ and $\mc{Y}$ satisfy
~\cite[(S1,2)]{Versal} and their deformations satisfy
~\cite[(4.1.i),(4.1.iii)]{Versal}.  If for both $e$ and $f$ there
exist generic extenders of algebraic deformations compatible with
\'etale extensions and with closed points, the same holds for $f\circ
e$.
\end{prop}

\begin{proof}
Let $(\widetilde{q}_{\mc{Y}},\widetilde{q}_{\mc{Y},0},b_{\mc{W}})$ be
a deformation situation over $f\circ e$ whose associated extension
pair is $B'\xrightarrow{q} B\xrightarrow{q_0} B_0$.  Assume the kernel
$M$ of $q$ is a module over $B_0/\mf{a}$ for a radical ideal $\mf{a}$.
By Lemma~\ref{lem-extend}, there exists a generic extender
$(A_0,N_{f\circ e})$ for the deformation situation over $f\circ e$.
The problem is to prove the generic extender is compatible with
\'etale extensions and closed points. The main issue is to find an
appropriate finite type, $\mf{a}$-generic localization of $A_0$.

First of all, by ~\cite[(4.1.iii)]{Versal}, there is a finite type,
$\mf{a}$-generic localization of $B_0$ so that $D_f(M)$ and $D_{e\circ
f}(M)$ are compatible with closed points.  Replace $B_0$ by this
localization and replace $B$ and $B'$ by the associated localizations.
Denote by $\widetilde{q}_{\mc{W},N_{f\circ e}}$ a lift to $\mc{W}$ of
the base change of $\widetilde{q}_{\mc{Y}}$ to $B'_{N_{f\circ e}}$.
Replace $B_0$ by $A_0$ and replace $B$ and $B'$ by the associated
localizations.  By hypothesis there exists a generic extender
$(A_0,N_f)$ for $f$ compatible with \'etale extensions and with closed
points.  Replace $B_0$ by $A_0$, etc.  Let
$b'_{\mc{X}}\xrightarrow{\widetilde{q}_{\mc{X}}} b_{\mc{X}}$ denote a
lift to $\mc{X}$ of the base change of $\widetilde{q}_{\mc{Y}}$ to
$B'_{N_f}$.  Because the base change of $\widetilde{q}_{\mc{Y}}$ to
$B'_{N_{f\circ e}}$ lifts to $e(\widetilde{q}_{\mc{W},N_{f\circ e}})$,
$N_{f\circ e}$ is a quotient of $N_f$.

``Calibrate'' the lifts $\widetilde{q}_{\mc{W},N_{f\circ e}}$ and
$\widetilde{q}_{\mc{X}}$ as follows.  Denote by
$\widetilde{q}_{\mc{X},N_{f\circ e}}$ the base change of
$\widetilde{q}_{\mc{X}}$ to $B'_{N_{f\circ e}}$.  The difference of
$e(\widetilde{q}_{\mc{W},N_{f\circ e}})$ and
$\widetilde{q}_{\mc{X},N_{f\circ e}}$ is an element $d$ of
$D_{f,\widetilde{q}_{\mc{W},0},b_{\mc{X},0}}(N_{f\circ e})$.  After a
further finite type, $\mf{a}$-generic localization, the surjection
$N_f\rightarrow N_{f\circ e}$ is split.  Therefore there exists an
element $d'$ of $D_{f,\widetilde{q}_{\mc{W},0},b_{\mc{X},0}}(N_f)$
mapping to $d$.  Replace $\widetilde{q}_{\mc{X}}$ by its translation
by $d$.  Then $e(\widetilde{q}_{\mc{W},N_{f\circ e}})$ equals
$\widetilde{q}_{\mc{X},N_{f\circ e}}$.

\begin{claim}  \label{claim-1}
It suffices to prove the proposition after replacing $M$ by the
quotient $N_f$, replacing $B'$ by $B'_{N_f}$, etc.
\end{claim}

For every \'etale extension $B'\rightarrow B'_{\text{\'et}}$ and every
projective, extending quotient $M\otimes_{B_0}
B_{\text{\'et},0}\rightarrow N_{\text{\'et}}$ for $f\circ e$ over
$B'_{\text{\'et}}$, also $N_{\text{\'et}}$ is extending for $f$ over
$B'_{\text{\'et}}$.  By hypothesis, $N_f$ is compatible with \'etale
extension.  Thus $M_{\text{\'et}}$ is a quotient of $N_f$.

Similarly, for every closed point $y$ of $\SP B_0$, and every
surjection $M\otimes_{B_0}k(y) \rightarrow N_y$, if the base change of
$\widetilde{q}_{\mc{Y}}$ to $B'_{N_y}$ lifts to $\mc{W}$, then it also
lifts to $\mc{X}$.  By hypothesis, $N_f$ is compatible with closed
points.  Thus $N_y$ is a quotient of $N_f$.  This proves Claim
~\ref{claim-1}.  Therefore, without loss of generality, replace $M$ by
the quotient $N_f$, replace $B'$ by $B'_{N_f}$, etc.

Here is the problem: the lift $\widetilde{q}_{\mc{X}}$ is not unique.
For every lift, by hypothesis, there is a finite type,
$\mf{a}$-generic localization and a generic extender over this
localization which is compatible with etale extensions and finite
points.  However, the finite type localization would appear to depend
on the choice of lift.  As there may be infinitely many lifts, it is a
priori possible the coproduct of these localizations is not finite
type.

Here is the remedy.  Let $D$ denote the module
$D_{f,\widetilde{q}_{\mc{Y},0},b_{\mc{X},0}}(B_0/\mf{a})$.  Localize
so that $D$ is a finite projective $B_0/\mf{a}$-module.  Denote the
module $\text{Hom}_{B_0}(D,B_0/\mf{a})$ by $D^\vee$.  Then
$D_{f,\widetilde{q}_{\mc{Y},0},b_{\mc{X},0}}(D^\vee)$ is canonically
isomorphic to $\text{Hom}_{B_0}(D,D)$.  Let $\beta$ denote the object
of $\mc{X}$ over $B_0 + D^\vee$ corresponding to $\text{Id}_D$.  Of
course $(B_0+D^\vee)\times_{B_0} B'$ equals $B'+D^\vee$.  By
~\cite[(S1)]{Versal}, there exists an extension $b_{\mc{X},D}'
\xrightarrow{\widetilde{q}_{\mc{X},D}} b_{\mc{X}}$ of
$\widetilde{q}_{\mc{X}}$ over $B'+D^\vee$ such that the image over
$B'$ is $\widetilde{q}_{\mc{X}}$ and the image over $B_0 + D^\vee$ is
$\beta$.

Here is the point: for every lift $b_{\mc{X},1}'
\xrightarrow{\widetilde{q}_{\mc{X},1}} b_{\mc{X}}$ of
$\widetilde{q}_{\mc{Y}}$, the difference of $b_{\mc{X},1}'$ and
$b_{\mc{X},1}$ is an element of
$D_{f,\widetilde{q}_{\mc{Y},0},b_{\mc{X},0}}(M)$, i.e., an element $d$
of $\text{Hom}_{B_0}(D^\vee,M)$.  There is a unique surjection of
$B'$-algebras, $B'+D^\vee \rightarrow B'$ whose restriction to
$D^\vee$ is $d:D^\vee \rightarrow M$.  The image of $b_{\mc{X},D}'$
equals $b_{\mc{X},1}'$.

By hypothesis, for the ideal $\mf{a}$ the deformation situation
$(\widetilde{q}_{\mc{X},D},\widetilde{q}_{\mc{X}},b_{\mc{W}})$ over
$e$ has a generic extender $(A_0,N_D)$ compatible with \'etale
extensions and closed points.  Let $D^\vee \rightarrow E$ be the fiber
coproduct of $M\oplus D^\vee \rightarrow N_D$ and the projection
$\text{pr}_{D^\vee}:M\oplus D^\vee \rightarrow D^\vee$.  After a
further finite type, $\mf{a}$-generic localization, $E$ is projective.

\begin{claim} \label{claim-2}
The quotient $M\oplus D^\vee \rightarrow N_D$ is the direct sum of the
quotient $M\rightarrow N_{f\circ e}$ and the quotient $D^\vee
\rightarrow E$.
\end{claim}

The claim holds if and only if it holds after a finite type,
$\mf{a}$-generic localization.

Let $K_D$ denote the kernel of the surjection $M\oplus D^\vee
\rightarrow N_D$.  Let $K_{f\circ e}$ denote the kernel of the
surjection $M\rightarrow N_{f\circ e}$.  The composition of
$\text{pr}_M:M\oplus D^\vee \rightarrow M$ and $M\rightarrow N_{f\circ
e}$ is extending.  Therefore, after a further $\mf{a}$-generic
localization, it factors through $N_D$, i.e., the image
$\text{pr}_M(K_D)$ is contained in $K_{f\circ e}$.  On the other hand,
the quotient $M/\text{pr}_M(K_D)$ is extending.  After a finite type,
$\mf{a}$-generic localization, it is projective.

Because $M \rightarrow M/\text{pr}_{M}(K_D)$ is extending and locally
free, after a further localization it is a quotient of $N_{f\circ e}$,
i.e., $K_{f\circ e}$ is contained in $\text{pr}_{M}(K_D)$.  Therefore
$\text{pr}_M(K_D)$ equals $K_{f\circ e}$.

Denote by $K'_D$ the kernel of the surjection $K_D\rightarrow
K_{f\circ e}$.  There is a projection
$\text{pr}_{K,D^\vee}:K'_D\rightarrow D^\vee$.  After a finite type,
$\mf{a}$-generic localization, the cokernel $E$ is projective.  The
composition $K_D \xrightarrow{\text{pr}_{D^\vee}} D^\vee \rightarrow
E$ factors through $K_D\rightarrow K_{f\circ e}$.  Denote this
morphism by $\delta:K_{f\circ e} \rightarrow E$.  The claim is
equivalent to the vanishing of $\delta$.

By way of contradiction, assume $\delta$ is not zero.  Then there
exists a surjection $E\rightarrow Q$ such that $K_{f\circ
e}\rightarrow Q$ is surjective and $Q$ is a projective
$B_0/\mf{a}$-module of length $1$.  Denote by $K_{f\circ e,Q}$ the
kernel of $K_{f\circ e}\rightarrow Q$.  Take the quotient of $M$ by
$K_{f\circ e,Q}$.  There is an induced isomorphism $K_{f\circ
e}/K_{f\circ e,Q} \rightarrow Q$.  Inverting this isomorphism and
composing with the inclusion gives a map of $B_0/\mf{a}$-modules,
$i:Q\rightarrow M/K_{f\circ e,Q}$.  Define $i':Q\rightarrow M$ to be
any $B_0/\mf{a}$-morphism lifting $i$.  Define $\phi:D^\vee
\rightarrow M$ to be the composition of $D^\vee \rightarrow Q$ and the
\emph{negative} of $i'$.  Define $(\text{Id}_M,\phi):M\oplus D^\vee
\rightarrow M$ to be the induced surjection.  The induced morphism
$(\text{Id}_M,\phi):K_D\rightarrow M$ annihilates $K'_D$.  So there is
an induced map from $K_{f\circ e}$ to $M$.  The restriction to
$K_{f\circ e,Q}$ is the inclusion.  So there is an induced map from
$K_{f\circ e}/K_{f\circ e,Q}$ to $M/K_{f\circ e,Q}$.  But, by
construction, this is the zero map.  Therefore
$(\text{Id}_M,\phi)(K_D)$ is strictly contained in $K_{f\circ e}$.
The quotient $M \rightarrow M/(\text{Id}_M,\phi)(K_D)$ is extending,
but does not factor through $M\rightarrow N_{f\circ e}$.  This
contradiction proves that $\delta$ is zero, i.e., it proves Claim
~\ref{claim-2}.

\begin{claim} \label{claim-3}
After this sequence of localizations, the generic extender
$M\rightarrow N_{f\circ e}$ for $f\circ e$ is compatible with \'etale
extensions and with closed points.
\end{claim}

\medskip\noindent
\textbf{\'Etale extensions.}  By hypothesis, the generic extender
$M\oplus D^\vee \rightarrow N_D$ is compatible with \'etale
extensions.  By Claim~\ref{claim-3}, this surjection is split.  It
follows that $M\rightarrow M_{f\circ e}$ is also compatible with
\'etale extensions.

\medskip\noindent
\textbf{Closed points.}  Let $B_0 \rightarrow k(y)$ be a closed point,
let $M\otimes_{B_0} k(y) \rightarrow N_y$ be a surjection, and let
$\widetilde{q}_{\mc{W},y}$ be a lift to $\mc{W}$ of the base change of
$\widetilde{q}_{\mc{Y}}$ to $B'_{N_y}$.  Because $D_f(M)$ is
compatible with closed points, there is a morphism $d:D^\vee
\rightarrow M$ so that $e(\widetilde{q}_{\mc{W},y})$ is the base
change of $\widetilde{q}_{\mc{X},D}$ associated to the composition of
$(\text{Id}_M,d):M\oplus D^\vee \rightarrow M$ and $M\rightarrow N_y$.
Because $M\oplus D^\vee \rightarrow N_D$ is compatible with closed
points, $N_y$ is a quotient of $N_D$, i.e., the kernel of
$M\rightarrow N_y$ contains $(\text{Id}_M,d)(K_D)$.  Because $K_D$ is
the direct sum $K_{f\circ e}\oplus K'_D$, in particular
$(\text{Id}_M,d)(K_D)$ contains $(\text{Id}_M,d)(K_{f\circ e}\oplus
0)$, i.e., $K_{f\circ e}$.  Thus $N_y$ factors through $M\rightarrow
N_{f\circ e}$.  Therefore $M\rightarrow N_{f\circ e}$ is compatible
with closed points.
\end{proof}

\subsection{Openness of versality} \label{ssect-thm}
Let $\mc{X}$ be a stack in groupoids over $(\text{Aff}/S)$ for the
\'etale topology.

\begin{defn} \label{defn-ov}
The stack $\mc{X}$ satisfies \emph{openness of versality} if for every
finite type $\OO_S$-algebra $R$ and every object $v$ of $\mc{X}(\SP
R)$, if $v$ is formally versal at $p$, then there is an open
neighborhood of $p$ in $\SP R$ on which $v$ is formally smooth.
Openness of versality is also \emph{compatible with \'etale
extensions} if for every \'etale morphism $e:\SP R^*\rightarrow \SP
R$, for every point $p^*$ of $\SP R^*$, if $v$ is formally versal at
$e(p^*)$, then $e^*v$ is formally smooth on an open neighborhood of
$p^*$ in $\SP R^*$.
\end{defn}

\begin{thm}~\cite[Proposition 4.3 and Theorem 4.4]{Versal} \label{thm-44}
Assume $\mc{X}$ is limit preserving, deformations of $\mc{X}$ satisfy
~\cite[(S1,2), (4.1)]{Versal}, every generic lift of an integral,
algebraic deformation situation over $\mc{X}$ has an integral lift,
and every algebraic deformation situation over $\mc{X}$ has a generic
extender compatible with \'etale extensions and with closed points,
cf. Definitions ~\ref{defn-lift} and ~\ref{defn-extend}.  Then
$\mc{X}$ satisfies openness of versality compatibly with \'etale
extensions.
\end{thm}

\begin{rmk} \label{rmk-44}
The hypotheses of Theorem~\ref{thm-44} do not include existence or
properties of an obstruction theory for $\mc{X}$.  In practice, an
obstruction theory is used to verify existence of integral lifts and
generic extenders, etc.
\end{rmk}

\begin{proof}
The entire proof will not be repeated.  We only indicate the
modifications necessary to replace existence and properties of
obstruction theory by the hypotheses above.

First of all, the proof of ~\cite[Proposition 4.3]{Versal} uses
obstructions only at the end of the proof.  In fact, obstructions are
used precisely to prove existence of generic extenders compatibly with
\'etale extensions.

Next, in the proof of ~\cite[Theorem 4.4]{Versal}, obstructions are
used to prove formal versality is stable under generization.  The
arguments on ~\cite[p. 176]{Versal} use obstructions precisely to
deduce existence of an integral lift $b'$ of the generic lift $a'$ of
the deformation situation $(B'\xrightarrow{q} B\xrightarrow{q_0} B_0,
b)$.  By ~\cite[Theorem 3.3]{Versal}, formal versality of $v$ at $p$
implies versality of $v$ at $p$.  Thus, there is a morphism from $R$
to the Henselization of the localization of $B'$ at $p$ inducing a
lift from $v$ to the ``Henselization'' of $b'$ at $p$, compatibly with
the map from $v$ to the ``Henselization'' of $b$ at $p$.  Because
$\mc{X}$ is limit preserving, this morphism factors through an \'etale
extension $B'\rightarrow B'_{\text{new}}$.  Replace the original
deformation situation by the base change by this \'etale extension.

Because the base change of $v$ by $R\rightarrow B'$ agrees with $b'$
after localizing at $p$, and by the same sort of arguments as above,
after replacing $B'$ by a subring $B'[t^{-1}M]$ of $A'$, the base
change of $v$ equals $b'$.  So the composition of $R\rightarrow B'$
with the localization $B'\rightarrow A'$ gives a morphism $v
\rightarrow a'$ making the diagram commute.  Therefore $v$ is formally
versal at the point $x$.

The final use of obstructions in the proof of ~\cite[Theorem
  4.4]{Versal} is at the bottom of p. 177 and the top of p. 178.  Here
obstructions are used to construct a generic extender $(A_0,N)$
compatibly with closed points.  Replace $B_0$ by $A_0$, and replace
$B$, $B'$ by the associated localizations.  By the hypothesis on
$\mc{S}$, the quotient $N$ has nonzero fiber $N\otimes_{B_0} k(y)$ for
a dense set of points.  Therefore the generic fiber of $N$ is nonzero.
Let $M\xrightarrow{\phi} B_0$ be any surjection factoring through
$M\rightarrow N$.  Denote by $B'\rightarrow B^*$ the associated
surjection.  This has the desired property.  The rest of the proof
goes through precisely as in ~\cite{Versal}; there is no further use
of obstructions.
\end{proof}

Let there be given a sequence of 1-morphisms of stacks in groupoids
over $(\text{Aff}/S)$,
$$
\mc{X}_n \xrightarrow{f_{n-1}^n} \mc{X}_{n-1} \xrightarrow{f_{n-2}^{n-1}} \dots
\xrightarrow{f_{1}^{2}} \mc{X}_1 \xrightarrow{f_{0}^{1}} (\text{Aff}/S).
$$ 
Assume each $\mc{X}_i$ satisfies ~\cite[(S1,2)]{Versal}.

\begin{cor} \label{cor-main}
For every $i=1,\dots,n$, assume $\mc{X}_i$ is limit preserving.  Also,
for every $i=1,\dots,n$, assume $f_{i-1}^i$ satisfies the relative
version of ~\cite[(S1,2)]{Versal} and the conditions of
Notation~\ref{notat-41a}.  Finally, for every $i=1,\dots,n$, assume
every generic lift of an integral, algebraic deformation situation
over $f_{i-1}^i$ has an integral lift, and every algebraic deformation
situation over $f_{i-1}^i$ has a generic extender compatible with
\'etale extensions and with closed points.  Then every $\mc{X}_i$
satisfies openness of versality compatibly with \'etale extensions.
\end{cor}

\begin{proof}
This follows immediately from Lemmas ~\ref{lem-sorites},
~\ref{lem-sorites2}, ~\ref{lem-sorites3} and Propositions
~\ref{prop-comp1}, ~\ref{prop-comp2}.
\end{proof}

\subsection{Artin's representability theorems} \label{ssect-rep}
Artin's approximation theorem has been generalized by Conrad and de
Jong, ~\cite{CdJ}.  This gives the following version of Artin's
representability theorem.

\begin{cor}~\cite[Corollary 5.2]{Versal}, ~\cite[Corollaire
    10.8]{LM-B}, ~\cite[Theorem 1.5]{CdJ}
  \label{cor-52}  
Assume $S$ is an excellent scheme.  Let $\mc{X}$ be a limit preserving
stack in groupoids over $(\text{Aff}/S)$ with the \'etale topology.
Then $\mc{X}$ is an algebraic stack if
\begin{enumerate}
\item[(1)]
$\mc{X}$ is relatively representable, i.e., the diagonal 1-morphism
  $\Delta:\mc{X} \rightarrow \mc{X}\times\mc{X}$ is representable by
  locally finitely presented algebraic spaces.
\item[(2)]
The Schlessinger-Rim criteria ~\cite[(S1,2)]{Versal} hold.
\item[(3)] 
If $\widehat{A}$ is a complete local $\OO_S$-algebra with residue
field of finite type over $S$, then $\mc{X}(\SP \widehat{A})
\rightarrow \varprojlim \mc{X}(\SP \widehat{A}/\mf{m}^n)$ has a dense
image.
\item[(4)]
The stack $\mc{X}$ satisfies openness of versality compatibly with
\'etale extensions.
\end{enumerate}
\end{cor}

Applying Theorem~\ref{thm-44}, this gives the following version.

\begin{thm} ~\cite[Theorem 5.3]{Versal}
  \label{thm-53}  
Assume $S$ is an excellent scheme.  Let $\mc{X}$ be a limit preserving
stack in groupoids over $(\text{Aff}/S)$ with the \'etale topology.
Then $\mc{X}$ is an algebraic stack if
\begin{enumerate}
\item[(1)]
The Schlessinger-Rim criteria ~\cite[(S1,2)]{Versal} hold.  Also, for
every algebraic element $a_0$ of $\mc{X}(\SP A_0)$ and every finite
$A_0$-module $M$, $\text{Aut}_{a_0}(A_0+M)$ is a finite $A_0$-module.
\item[(2)]
For any complete local $\OO_S$-algebra $\widehat{A}$ with residue
field of finite type over $S$, the canonical map
$$
\mc{X}(\SP \widehat{A})\rightarrow \varprojlim \mc{S}(\SP
\widehat{A}/\mf{m}^n) 
$$
is faithful, and has a dense image, i.e., the projection to
$\mc{X}(\SP \widehat{A}/\mf{m}^n)$ is essentially surjective for every
$n$.
\item[(3)]
$D$ and $\text{Aut}_{a_0}(A_0+M)$ satisfy ~\cite[(4.1)]{Versal}, 
every generic lift of an integral, algebraic deformation situation
over $\mc{X}$ has an integral lift, and every algebraic deformation
situation over $\mc{X}$ has a generic extender compatible with \'etale
extensions and with closed points.
\item[(4)]
If the object $a_0$ of $\mc{X}(\SP A_0)$ is algebraic, and if $\phi$
is an automorphism of $a_0$ inducing the identity in $\mc{X}(\SP
k(y))$ for a dense set of finite type points $y$ of $\SP A_0$, then
$\phi$ equals $\text{Id}_{a_0}$ on a non-empty open subset of $\SP
A_0$.
\end{enumerate}
\end{thm}

\subsection{Relative obstruction theory} \label{ssect-obs}
Relative obstruction theories are studied in ~\cite[Appendix A]{OK3}.
Let $f:\mc{X}\rightarrow \mc{Y}$ be a 1-morphism between stacks in
groupoids over $(\text{Aff}/S)$ both satisfying Schlessinger's
conditions ~\cite[(S1,2)]{Versal}.

\begin{defn}~\cite[Def. A.10]{OK3} \label{defn-obs}
A \emph{relative obstruction theory for $f$} consists of the
  following:
\begin{enumerate}
\item[(i)]
an assignment to each reduced infinitesimal extension
$(\widetilde{q}_0,a)$ over $f$ of an $A_0$-linear functor,
$$
O_{\widetilde{q}_0,a}: (\text{finite } A_0-\text{modules}) \rightarrow
(\text{finite } A_0-\text{modules}), 
$$
\item[(ii)]
an assignment to each morphism $(\widetilde{u},\widetilde{u}_0,\phi)$
of infinitesimal extensions over $f$ of a natural transformation of
$A_0$-linear functors
$$
\theta_{\widetilde{u},\widetilde{u}_0,\phi}:
O_{\widetilde{q}_{A,0},a}(M)\otimes_{A_0} B_0 \Rightarrow 
O_{\widetilde{q}_{B,0},b}(M\otimes_{A_0} B_0),
$$
\item[(iii)]
and an assignment to each deformation situation
$(\widetilde{q},\widetilde{q}_0,a)$ over $f$ of an element
$o_{\widetilde{q},\widetilde{q}_0,a}$ of $O_{\widetilde{q},a}(M)$
\end{enumerate}
satisfying the following axioms.
\begin{enumerate}
\item[(i)]
For every deformation situation $(\widetilde{q},\widetilde{q}_0,a)$
over $f$, $o_{\widetilde{q},\widetilde{q}_0,a}$ equals $0$ if and only
if there exists a morphism $a'\rightarrow a$ in $\mc{X}$ mapping to
$\widetilde{q}$.
\item[(ii)]
For every morphism
$(\widetilde{u}',\widetilde{u},\widetilde{u}_0,\phi)$ of extension
pairs over $f$, $\theta_{\widetilde{u},\widetilde{u}_0,\phi}(
o_{\widetilde{q}_A,\widetilde{q}_{A,0},a})$ equals
$o_{\widetilde{q}_B,\widetilde{q}_{B,0},b}$.
\end{enumerate}
\end{defn}

\begin{notat}~\cite[(4.1)]{Versal}, ~\cite[\S A.11]{OK3} \label{notat-41}
As in Notation~\ref{notat-41a}, there are relative analogues of the
conditions on obstruction theory.
\begin{enumerate}
\item[(4.1.i)]
The functor $O_{f}$ is compatible with \'etale localization: For every
morphism $(\widetilde{u},\widetilde{u}_0,\phi)$ of infinitesimal
extensions over $f$, if $A_0$ is a finite-type $\OO_S$-algebra and $u$
is \'etale, then the following associated natural transformation of
functors is an isomorphism,
$$
O_{f,\widetilde{q}_{A,0},a_0}(M)\otimes_{A_0} B_0 \rightarrow
O_{f,\widetilde{q}_{B,0},b_0}(M\otimes_{A_0} B_0).
$$
\item[(4.1.iii)]
For every infinitesimal extension $(\widetilde{q}_0,a)$ over $f$ with
$A$ a finite type $\OO_S$-algebra, there is an open dense set of
points of finite type $p\in \SP A_0$ so that the following map is
injective
$$
O_{\widetilde{q}_0,a}(M)\otimes_{A_0} k(p) 
\rightarrow O_{\widetilde{q}_0,a}(M\otimes_{A_0} k(p)).
$$
\end{enumerate}
\end{notat}

\subsection{Obstructions and lifting properties} \label{ssect-verify}
Properties of obstructions imply the hypotheses of Theorem
~\ref{thm-44}.
 
\begin{lem}~\cite[Lemma 4.6]{Versal} \label{lem-46a}
Let $f:\mc{X}\rightarrow \mc{Y}$ be a 1-morphism of stacks in
groupoids over $(\text{Aff}/S)$.  Assume both satisfy
~\cite[(S1,2)]{Versal}.  Let $O_f$ be a relative obstruction theory
for $f$.  If $\mc{X}$ is limit preserving, if $D_f$ satisfies (4.1.i)
of Notation~\ref{notat-41a} and if $O_f$ satisfies (4.1.i) of
Notation~\ref{notat-41}, then every generic lift of an integral,
algebraic deformation situation over $f$ has an integral lift.
\end{lem}

\begin{proof}
Denote by $M$ the kernel of $\widetilde{q}_B$.  The obstruction to
finding a morphism $r_B$ over $\widetilde{q}_B$ is an element $o$ of
the obstruction group $O=O_{f,\widetilde{q}_{B,0},b_{\mc{X}}}(M)$.  By
hypothesis, there is a lift after localizing to the total rings of
fractions.  These localizations are limits of finite type
localizations.  Because $\mc{X}$ is limit preserving, there exists a
lift after a finite type localization.  Thus, the image of $o$ in the
obstruction group of this finite type localization is zero.  Because a
finite type localization is an \'etale morphism, (4.1.i) of
Notation~\ref{notat-41} implies the obstruction group of the
localization is the localization of the obstruction group.  Therefore,
the image of $o$ in $O\otimes_{B_0} A_0$ is zero.  So there exists a
nonzero element $t$ of $B_0$ such that $t\cdot o$ is zero in $O$.

Let $B'_{\text{new}}$ be the subring of $A'$ generated by $B'$ and by
$t^{-1}M$ inside $M\otimes_{B_0} A_0$, i.e., the kernel of $q_A$.
Denote by $q_{B,\text{new}}$ the restriction of $q_A$ to
$B'_{\text{new}}$ considered as a morphism with target $B$.  The
infinitesimal extension $B'_{\text{new}}
\xrightarrow{q_{B,\text{new}}} B$ gives an deformation situation
$(q_{B,\text{new}},q_{B,0})$ which is the image of an obvious morphism
from $(q_B,q_{B,0})$.  The base change of
$(\widetilde{q}_B,\widetilde{q}_{B,0},b_{\mc{X}})$ by this morphism is
an deformation situation
$(\widetilde{q}_{B,\text{new}},\widetilde{q}_{B,0},b_{\mc{X}})$ over
$f$.  By Axiom (ii) of Definition ~\ref{defn-obs}, the obstruction
$o_{\text{new}}$ of this deformation situation is the image of $o$ in
$O_{f,\widetilde{q}_{B,0},b_{\mc{X}}}(t^{-1}M)$.  There is an
isomorphism $t^{-1}M \rightarrow M$ sending $t^{-1}m$ to $m$.  This
defines an isomorphism $O_{f,\widetilde{q}_{B,0},b_{\mc{X}}}(t^{-1}M)
\rightarrow O_{f,\widetilde{q}_{B,0},b_{mc{X}}}(M)$, and the image of
$o_{\text{new}}$ is exactly $t\cdot 0$, which is zero.  Therefore
$o_{\text{new}}$ is zero.  By Axiom (i) of Definition ~\ref{defn-obs},
there is a morphism $r^*_B$ of $\mc{X}$ over
$\widetilde{q}_{B,\text{new}}$.

Replace $B'$ by $B'_{\text{new}}$, etc.  By construction, $r^*_B$
satisfies Axioms (i) and (ii) of an integral lift.  The base change of
$r^*_B$ over $A'\rightarrow A$ is a morphism $r^*_A$ over
$\widetilde{q}_A$.  By ~\cite[(S1,2)]{Versal}, there exists an element
$\delta$ of $D_{f,a_{\mc{X},0}}(M\otimes_{B_0}A_0)$ such that
$\delta\cdot r^*_A$ equals $r_A$.  By the same type of argument as
above, using (4.1.i) of Notation~\ref{notat-41a} for $D_f$ and that
$\mc{X}$ is limit preserving, there exists a nonzero element $t$ of
$A_0$ and an element $d$ of $D_{f,b_{\mc{X},0}}(M)$ such that $t\delta
= d$.  Thus, after replacing $B'$ by an extension $B'\rightarrow
B'_{\text{new}}$ as above and replacing $r^*_B$ by its base change
over $B'_{\text{new}}$, there exists an element $d$ of
$D_{f,b_{\mc{X},0}}(M)$ such that $r_A$ equals $d\cdot r^*_A$.  Define
$r_B$ to be $d\cdot r^*_B$.  This is an integral lift of $r_A$.
\end{proof}

\begin{lem}~\cite[Proof of Theorem 4.4]{Versal} \label{lem-surj}
Let $f:\mc{X}\rightarrow \mc{Y}$ be a 1-morphism of stacks in
groupoids over $(\text{Aff}/S)$.  Assume both satisfy
~\cite[(S1,2)]{Versal}.  Let $O_f$ be a relative obstruction theory
for $f$.  If $O_f$ satisfies (4.1.i) and (4.1.iii) of
Notation~\ref{notat-41}, then every algebraic deformation situation
over $f$ has a generic extender compatible with \'etale extension and
with closed points.
\end{lem}

\begin{proof}
Because of (4.1.i), we may replace $B_0$ by any finite type
localization (and replace $B$ and $B'$ by the associated
localizations).  In particular, assume $M$ is a finite projective
$B_0$-module.  Similarly, assume $P =
O_{f,\widetilde{q}_0,b_{\mc{X}}}(B_0)$ is a finite projective
$B_0$-module.

By $B_0$-linearity of $O_f$, $O_{f,\widetilde{q}_0,b_{\mc{X}}}(M)$ is
canonically isomorphic to $M\otimes_{B_0}P$.  Moreover, for every
surjection $M\rightarrow M'$, the map $M\otimes_{B_0} P \rightarrow
O_{b,\widetilde{q}_0,b_{\mc{X}}}(M')$ factors through the surjection
$M\otimes_{B_0} P \rightarrow M'\otimes_{B_0} P$.  By (4.1.iii), after
replacing $B_0$ by a further finite type localization, the induced map
$M'\otimes_{B_0} P \rightarrow O_{f,\widetilde{q}_0,b_{\mc{X}}}(M')$
is injective whenever $M'$ is a direct sum of copies of $k(y)$ for $y$
a closed point.

The obstruction class $o$ associated to
$(\widetilde{q},\widetilde{q}_0,b_{\mc{X}})$ gives a $B_0$-linear map,
$B_0\rightarrow M\otimes_{B_0} P$.  By adjunction, this is equivalent
to a $B_0$-linear map, $\text{Hom}_{B_0}(P,B_0)\rightarrow M$.  Define
$M\rightarrow N$ to be the cokernel of this map.  For every surjection
$M\rightarrow N'$, the image of $o$ in $N'\otimes_{B_0} P$ is zero if
and only if $N'$ is a quotient of $N$.  Because the image of $o$ in
$N\otimes_{B_0} P$ is zero, the image is also zero in
$O_{f,\widetilde{q}_0,b_{\mc{X}}}(N)$.  Therefore there exists a
morphism $b_{\mc{X},N}\rightarrow b_{\mc{X}}$ in $\mc{X}$ mapping to
$\widetilde{q}_N$.  So $N$ is a generic extender.

For a quotient $M'$ of $M\otimes_{B_0}k(y)$.  the induced map
$M'\otimes_{B_0} P \rightarrow O_{f,\widetilde{q}_0,b_{\mc{X}}}(M')$
is injective.  Therefore the image of $o$ in
$O_{f,\widetilde{q}_0,b_{\mc{X}}}(M')$ is zero if and only if $M'$ is
a quotient of $N$.  Thus $N$ is compatible with closed points.

Finally, compatibility with \'etale extensions follows from (4.1.i).
\end{proof}

\begin{cor} \label{cor-algobs}
Assume $\mc{Y}$ is limit preserving and satisfies
~\cite[(S1,2)]{Versal} and $f$ is representable by limit preserving
algebraic stacks.  Then $\mc{X}$ is limit preserving and satisfies
~\cite[(S1,2)]{Versal}.  Also every generic lift of an integral,
algebraic deformation situation over $f$ has an integral lift, and
every algebraic deformation situation over $f$ has a generic extender
compatible with \'etale extension and with closed points.  Finally, if
$\mc{Y}$ is relatively representable, then $\mc{X}$ is relatively
representable.
\end{cor}

\begin{proof}
Let $A$ equal $\lim A_i$ and let $a_\mc{X}$ be an object of
$\mc{X}(\SP A)$.  Denote by $a_\mc{Y}$ the image $f(a_{\mc{X}})$.
Because $\mc{Y}$ is limit preserving, there exists an $i$ and an
object $a_{\mc{Y},i}$ of $\mc{Y}(\SP A_i)$ whose base change is
isomorphic to $a_{\mc{Y}}$.  Form the 2-fibered product,
$$
\mc{X}_{a_\mc{Y},i} =
\mc{X}\times_{f,\mc{Y},a_{\mc{Y},i}} \SP A_i.
$$
By hypothesis, this is a limit preserving.  And $a_{\mc{X}}$ gives an
object of $\mc{X}_{a_\mc{Y},i}(\SP A)$.  Thus there exists a $j$ and
an object $a_{\mc{X},j}$ of $\mc{X}_{a_\mc{Y},i}(\SP A_j)$ whose base
change is isomorphic to $a_{\mc{X}}$.  Therefore $\mc{X}$ is limit
preserving.  By a similar argument, $\mc{X}$ satisfies
~\cite[(S1,2)]{Versal}.

Given an algebraic deformation situation over $f$,
$(\widetilde{q},\widetilde{q}_0,b_{\mc{X}})$, the 2-fibered product
$$
\mc{X}_{b_\mc{Y}'}= \mc{X}\times_{f,\mc{Y},b_{\mc{Y}}'} \SP B'
$$ 
is a limit preserving algebraic stack.  Every base change of this
deformation situation over $f$ is the projection of a deformation
situation over $\mc{X}_{b_\mc{Y}'}$.  Thus it suffices to prove
algebraic deformation situations over $\mc{X}_{b_\mc{Y}'}$ have
integral lifts and generic extenders, etc.

By ~\cite[p. 182]{Versal}, ~\cite[Remark 1.7]{ODef},
$\mc{X}_{b_\mc{Y}'}$ has an obstruction theory such that
~\cite[(S1,2), (4.1)]{Versal} are satisfied for automorphisms,
deformations and obstructions.  By Lemmas ~\ref{lem-46a} and
~\ref{lem-surj}, algebraic deformation situations over
$\mc{X}_{b_\mc{Y}'}$ have integral lifts, generic extenders, etc.

Finally, if $\mc{Y}$ is relatively representable, then
$\mc{X}\times_{\mc{Y}}\mc{X} \rightarrow \mc{X}\times \mc{X}$ is
representable.  Since $f$ is relatively representable,
$\mc{X}\rightarrow \mc{X}\times_{\mc{Y}}\mc{X}$ is representable.
Therefore the composition $\mc{X}\rightarrow \mc{X}\times \mc{X}$ is
representable, i.e., $\mc{X}$ is relatively representable.
\end{proof}

\section{The stack of algebraic spaces} \label{sec-alg}

Denote by $\mc{X}$ the category whose objects are all pairs $(U,X)$ of
an affine scheme $U$ and a $U$-algebraic space $X$, and whose
morphisms are all pairs,
$$
(f_U,f_X):(U',X')\rightarrow (U'',X''),
$$
of a morphism $f_U:U'\rightarrow U''$ of affine schemes and an
isomorphism of $U'$-algebraic spaces, $f_X:X'\rightarrow
U'\times_{U''} X''$.  The identity morphisms and the composition law
are the obvious ones.  There is a functor $\mc{X}\rightarrow
(\text{Aff}/\SP \ZZ)$.  This functor is a stack for the \'etale
topology, and even the fpqc topology on $(\text{Aff})$,
cf. ~\cite[(3.4.6), (9.4)]{LM-B}.

\begin{claim} \label{claim-1a}
The diagonal 1-morphism $\Delta:\mc{X}\rightarrow \mc{X}\times \mc{X}$
is not representable.
\end{claim}

\begin{proof}
Let $k$ be any field.  Consider the 1-morphism $\zeta:\SP(k)
\rightarrow \mc{X}$ associated to the object $(\SP(k),\AAA^1_k)$.
There is an associated 1-morphism $(\zeta,\zeta):\SP(k) \rightarrow
\mc{X} \times \mc{X}$.  The 2-fibered product of $(\zeta,\zeta)$ and
$\Delta$ is the stack associated to the contravariant functor $I$ on
$(\text{Aff}/\SP(k))$ associating to every $k$-algebra $A$ the set
$$
I(A) := \text{Isom}_A(A[t], A[t]).
$$ 
The claim is that $I$ is not an algebraic space.
 
Consider the subset,
$$
T = \{ \phi \in I(k[\epsilon]/\epsilon^2) | \phi \equiv \text{Id}
(\text{mod} \epsilon) \}.
$$
By direct computation, there is an isomorphism from the $k$-vector
space $k[t]$ to $T$ via $f \mapsto (t\mapsto t+\epsilon f)$.  The
dimension of $k[t]$ as a $k$-vector space is countably infinite.
Therefore the dimension of $T$ as a $k$-vector space is countably
infinite.

On the other hand, for every algebraic space $J$ over $k$, and every
$k$-point $p$ of $J$, the set,
$$
T_{J,p} = \{ \phi \in J(k[\epsilon]/\epsilon^2) | \phi \equiv p
(\text{mod} \epsilon)\},
$$
is isomorphic to the $k$-vector space,
$$
\text{Hom}_k(\Omega,k),
$$
where $\Omega$ is the $k$-vector space
$\Omega_{J/k}/\mf{m}_p\Omega_{J/k}$.  If $T'$ is infinite-dimensional,
then $\Omega$ is also infinite-dimensional.  For an
infinite-dimensional $k$-vector space $\Omega$, the dimension of
$T'=\text{Hom}_k(\Omega,k)$ is uncountable.  Since the dimension of
$T$ is countably infinite, $I$ is not an algebraic space over $k$.
\end{proof}

Denote by $\mc{X}_\pflpf$ the category whose objects are all pairs
$(U,X)$ of an affine scheme $U$ and a proper, flat, locally finitely
presented $U$-algebraic space.  This is a full subcategory of
$\mc{X}$.  The restriction of $F$ is a functor $F:\mc{X}_\pflpf
\rightarrow (\text{Aff})$.  As with $\mc{X}$, $\mc{X}_\pflpf$ is a
stack for the \'etale and fpqc topology.

\begin{prop} \label{prop-pf}
The stack $\mc{X}_\pflpf$ is limit preserving.
\end{prop}

\begin{proof}
This follows from results in ~\cite[\S 8]{EGA4}.  
\end{proof}

\begin{prop} \label{prop-diag}
The diagonal morphism $\Delta:\mc{X}_\pflpf \rightarrow
\mc{X}_\pflpf\times \mc{X}_\pflpf$ is representable, separated and
locally finitely presented.
\end{prop}

\begin{proof}
This follows from ~\cite[Theorem 6.1]{Artin}.  
\end{proof}

Denote by $\pi:\mc{V}_\pflpf \rightarrow \mc{X}_\pflpf$ the universal
1-morphism representable by proper, flat, locally finitely presented
algebraic spaces.

\begin{prop} \label{prop-Artin}
Deformations and automorphisms of $\mc{X}_\pflpf$ satisfy
~\cite[(S1,2), (4.1)]{Versal}.  There is an obstruction theory
$\mc{O}$ for $\mc{X}_\pflpf$ satisfying ~\cite[(4.1)]{Versal}.
\end{prop}

\begin{proof}
The existence of an obstruction theory satisfying ~\cite[(S1)]{Versal}
follows from ~\cite[Proposition III.2.1.2.3]{Illusie1} using the
relative cotangent complex of $\pi$.  The
condition~\cite[(S2)]{Versal} for deformations and the analogous
condition for automorphisms follow from coherence of the cohomology
sheaves of the cotangent complex ~\cite[Corollaire II.2.3.7]{Illusie1}
together with the finiteness theorem, ~\cite[Th\'eor\`eme
III.2.2]{SGA6}.

Condition ~\cite[(4.1.i)]{Versal} follows from
~\cite[II.1.2.3.5]{Illusie1} and standard results about cohomology and
flat base change.  Condition ~\cite[(4.1.ii)]{Versal} follows from the
theorem on formal functions, ~\cite[Theorem V.3.1]{Kn}.  Condition
~\cite[(4.1.iii)]{Versal} follows from generic flatness,
~\cite[Th\'eor\`eme 6.9.1]{EGA4}, and the semicontinuity theorem
~\cite[Th\'eor\`eme 7.7.5]{EGA3}.
\end{proof}

However $\mc{X}_\pflpf$ does not satisfy Axiom 3 of ~\cite[Corollary
  5.2]{Versal}.

\begin{claim} \label{claim-neff}
There does not exist a pair $(Q,f)$ of an algebraic space $Q$ and a
representable, faithful, smooth 1-morphism $f:Q\rightarrow
\mc{X}_\pflpf$.  In fact, $\mc{X}_\pflpf$ does not satisfy Axiom 3 of
~\cite[Corollary 5.2]{Versal}.
\end{claim}

\begin{proof}
If $\mc{X}_\pflpf$ did satisfy Axiom 3, then every proper algebraic
space over a field would have an effective, formal, versal
deformation.  By ~\cite[Theorem 1.6]{Artin}, every proper algebraic
space over a field does have a formal, versal deformation.  However,
it is not always effective.  The following example of a projective,
smooth variety in characteristic $0$ with no effective, formal, versal
deformation is well-known.

Let $k$ be an uncountable, characteristic $0$, algebraically closed
field.  Let $X$ be a smooth anticanonical divisor in $\PP^1\times
\PP^2$.  This is a K3 surface together with an elliptic fibration
$\text{pr}_{\PP^1}:X\rightarrow \PP^1$.

Because the Schlessinger-Rim criteria ~\cite[(S1,2)]{Versal} hold for
$\mc{X}$, there is a complete, local $k$-algebra $R$ and a formal,
versal deformation $(X_n)_{n\geq 0}$ of $X$ over $R$.  Because K3
surfaces are unobstructed, $R$ is formally smooth, i.e., $R$ is a
power series ring.  Also there is a canonical isomorphism,
$$
\text{Hom}_k(\mathfrak{m}/\mathfrak{m}^2,k) \cong H^1(X,T_X).
$$
For every invertible sheaf $\mc{L}$ on $X$, there is an associated
first Chern class $C_1(\mc{L})$ in $H^1(X,\Omega_X)$.  This defines an
injective group homomorphism,
$$
C_1: N^1(X) \rightarrow H^1(X,\Omega_X),
$$
where $N^1(X)$ is the group of numerical equivalence classes of
invertible sheaves.

There is a cup-product pairing,
$$
H^1(X,\Omega_X) \times H^1(X,T_X) \rightarrow H^2(X,\OO_X).
$$
Because $X$ is a K3 surface, or equivalently by adjunction for the
inclusion of $X$ in $\PP^1 \times \PP^2$, there exists an isomorphism
of the dualizing sheaf $\omega_X$ with $\OO_X$.  Using this
isomorphism, the cup-product pairing above is equivalent to the
pairing for Serre duality, which is a perfect pairing.

Let $\mc{L}$ be an invertible sheaf on $X$ and let
$\theta:\mf{m}/\mf{m}^2 \rightarrow k$ be an element in $H^1(X,T_X)$,
considered as a $k$-algebra homomorphism $R/\mf{m}^2 \rightarrow
k[\epsilon]/\epsilon^2$.  Denote by $X_\theta$ the Abelian scheme over
$k[\epsilon]/\epsilon^2$ given by $\theta(X_2)$.  The invertible sheaf
$\mc{L}$ is the restriction of an invertible sheaf on $X_\theta$ only
if $C_1(\mc{L})\cup \theta$ is $0$ in $H^2(X,T_X)$.  If $\mc{L}$ is
not numerically trivial, then $C_1(\mc{L})$ is nonzero.  Then, because
the cup product pairing is nondegenerate, set of $\theta$ for which
$\mc{L}$ extends over $X_\theta$ is a proper $k$-subspace $V_{\mc{L}}$
of $H^1(X,T_X)$.

By the theorem of the base, the group of $N^1(X)$ of numerical
equivalence classes of invertible sheaves is a finitely generated
Abelian group.  Therefore, the set of subspaces $V_{\mc{L}}$ arising
from invertible sheaves as above is countable.  Because $k$ is
uncountable, there exists an element $\theta$ of $H^1(X,T_X)$
contained in none of the countably many proper subspaces $V_{\mc{L}}$.
Therefore every invertible sheaf on $X_\theta$ is numerically trivial.

Because $R$ is a power series ring, there is a local $k$-algebra
homomorphism $t:R\rightarrow k\Sem{\epsilon}$ whose associated map
$\mf{m}/\mf{m}^2 \rightarrow k$ is $\theta$.  If $(R/\mf{m}^n,X_n)$
comes from a proper algebraic space $X_R$ over $\SP(R)$, then $t(X_R)$
is a proper, flat algebraic space over $\SP k\Sem{\epsilon}$.  By
~\cite[Corollary 5.20]{Kn}, there is a dense open subspace $U$ of
$t(X_R)$ which is an affine scheme.  Denote by $D'$ the complement of
$U$ in $t(X_R)$.  Denote by $D$ the closure in $t(X_R)$ of the generic
fiber of $D'\rightarrow \SP k\Sem{\epsilon}$.  Since the closed fiber
of $t(X_R)$ is smooth, $t(X_R)$ is smooth over $\SP k\Sem{\epsilon}$.
Also, $t(X_R)$ is separated.  Therefore, essentially by
~\cite[Th\'eor\`eme 5.10.5]{EGA3}, $D$ is a Cartier divisor in $X_R$.
Moreover, it is flat over $\SP k\Sem{\epsilon}$.  Denote by $\mc{L}$
the associated invertible sheaf.  By the paragraph above, the
restriction of $\mc{L}$ to $X$ is numerically trivial.

Because the generic fiber of $t(X_R)$ is a proper,
positive-dimensional algebraic space, $U$ is not all of the generic
fiber, i.e., the generic fiber of $D'$ is not empty.  Because $D$ is
proper over $\SP k\Sem{\epsilon}$ and its image contains the generic
point, the closed fiber of $D$ is not empty.  Because $D$ is flat over
$\SP k\Sem{\epsilon}$, the closed fiber of $D$ is not all of $X$.
Thus it is a nonempty curve in $X$.  But no nonempty curve in $X$ is
numerically trivial: it has positive intersection number with either
the pullback of $\OO_{\PP^1}(1)$ or the pullback of $\OO_{\PP^2}(1)$.
This contradicts that the restriction of $\mc{L}$ to $X$ is
numerically trivial.  The contradiction proves there exists no proper
algebraic space $X_R$ giving the formal, versal deformation of $X$.
\end{proof}

\begin{rmk} \label{rmk-neff}
The counterexample also proves there is no effective, formal, versal
deformation of the pair $(X,\text{pr}_{\PP^1}^*\OO_{\PP^1}(1))$.  If
there were, then the argument above would prove there is a proper,
flat algebraic space $t(X_R)$ over $\SP k\Sem{\epsilon}$ and an
invertible sheaf $\mc{M}$ on $t(X_R)$ restricting to
$\text{pr}_{\PP^1}^* \OO_{\PP^1}(1)$ such that the restriction to $X$
of every invertible sheaf on $t(X_R)$ is numerically equivalent to
$\text{pr}_{\PP^1}^* \OO_{\PP^1}(d)$ for some $d$.  Because
$h^1(X,\text{pr}_{\PP^1}^*\OO_{\PP^1}(1))$ is zero, every global
section of $\text{pr}_{\PP^1}^*\OO_{\PP^1}(1)$ is the restriction of a
global section of $\mc{M}$.  Because
$\text{pr}_{\PP^1}^*\OO_{\PP^1}(1)$ is generated by global sections,
there exists a global section of $\mc{M}$ that does not vanish
identically on $D$.  The generic fiber of the zero locus of this
section is a nonempty curve not contained in $D$. Because it is
proper, it is also not contained in $U$.  Thus it has positive
intersection number with $D$.  On the other hand,
$C_1(\text{pr}_{\PP^1}^*\OO_{\PP^1}(1))^2$ is zero.  This proves the
restriction of $\mc{L}$ is not numerically equivalent to
$\text{pr}_{\PP^1}^*\OO_{\PP^1}(1)$.  The contradiction proves there
is no effective, versal, deformation of
$(X,\text{pr}_{\PP^1}^*\OO_{\PP^1}(1))$.
\end{rmk}

The following is also noteworthy.

\begin{claim} \label{claim-nqc}
The diagonal morphism is not quasi-compact.
\end{claim}

\begin{proof}
Let $k$ be an algebraically closed field and let $(E,0)$ be an
elliptic curve over $\SP(k)$ whose automorphism group is $\ZZ/2\ZZ$.
Let $X=E\times E$.  For the object $(\SP(k),X)$, the fiber in the
diagonal is the scheme $I=\text{Isom}_k(X,X)$.  The scheme $I$
surjects smoothly to $X$ by associating to each map $f$ the image
$f(0)$.  The kernel is isomorphic to the discrete, non-quasi-compact
$k$-group scheme $\text{GL}_2(\ZZ)$ by associating to $f$ the unique
matrix
$$
\lt( 
\begin{array}{cc}
a & b \\
c & d 
\end{array} \rt),
$$
in $\text{GL}_2(\ZZ)$ such that $f(x,y) = (ax+by,cx+dy)$.  
\end{proof}

\section{Variants of the stack of algebraic spaces}
\label{sec-variants}

As above, denote by $\pi:\mc{V}_\pflpf \rightarrow \mc{X}_\pflpf$ the
universal 1-morphism representable by proper, flat, locally finitely
presented algebraic spaces.  Denote by $\text{Coh}_\pi$ the category
whose objects consist of triples $(U,X,E)$ of an object $(U,X)$ of
$\mc{X}_\pflpf$ and a locally finitely presented $\OO_X$-module $E$
which is $\OO_U$-flat, cf. ~\cite[2.4.4]{LM-B}.  Morphisms in
$\text{Coh}_\pi$ are data $(f_U,f_X,f_E)$ of a morphism $(f_U,f_X)$ in
$\mc{X}_\pflpf$ and an isomorphism $f_E$ from $E'$ to the pullback of
$E''$.  There is a forgetful functor $G:\text{Coh}_\pi \rightarrow
\mc{X}_\pflpf$.

\begin{prop}~\cite[Th\'eor\`eme 4.6.2.1]{LM-B} \label{prop-Pic}
The category $\text{Coh}_\pi$ is a stack for the \'etale and fpqc
topologies on $(\text{Aff})$.  Moreover, the functor $G$ is
representable by limit preserving algebraic stacks with quasi-compact,
separated diagonal.
\end{prop}

\begin{proof}
By the proof of the first half of ~\cite[Th\'eor\`eme 4.6.2.1]{LM-B},
the diagonal of $\text{Coh}_{X/U}$ is representable by separated,
finitely presented algebraic spaces.  Note the first half of the proof
uses only properness of $X/U$ and does not use cohomological flatness
in dimension $0$.  To be absolutely precise, one has to generalize
~\cite[Corollaire 7.7.8]{EGA3} to the case when $f:X\rightarrow Y$ is
a proper morphism of algebraic spaces.  Using ~\cite{Versal}, it is
clear how to do this (perhaps ~\cite[Remarques 7.7.9.iii]{EGA3}
anticipates this generalization).

The proof that $\text{Coh}_{X/U}$ is a locally finitely presented
Artin stack uses \cite[Theorem 5.3]{Versal}.  Existence of an
obstruction theory for $\text{Coh}_{X/U}$ such that automorphisms,
deformations and obstructions satisfy ~\cite[(S1,2), (4.1)]{Versal} is
just as in Proposition~\ref{prop-Artin}, replacing \cite[Proposition
III.2.1.2.3]{Illusie1} by \cite[Proposition IV.3.1.5]{Illusie1}.
Finally, compatibility with completions follows from the Grothendieck
existence theorem, ~\cite[Theorem V.6.3]{Kn}.
\end{proof} 

Denote by $\mc{X}_\text{pol}$ the full subcategory of $\text{Coh}_\pi$
of triples $(U,X,L)$ where $L$ is an invertible sheaf on $X$ ample
relative to $U$.

\begin{prop} \label{prop-pol}
The category $\mc{X}_\text{pol}$ is a limit preserving algebraic stack
with quasi-compact, separated diagonal.
\end{prop}

\begin{proof}
By ~\cite[Th\'eor\`eme 4.7.1]{EGA3}, the inclusion functor
$\mc{X}_\text{pol} \rightarrow \text{Coh}_\pi$ is representable by
open immersions.  Therefore $\mc{X}_\text{pol} \rightarrow
\mc{X}_\pflpf$ is representable by limit preserving algebraic stacks.
By Proposition ~\ref{prop-comp1}, Proposition ~\ref{prop-comp2},
Corollary ~\ref{cor-algobs}, Proposition~\ref{prop-pf} and
Proposition~\ref{prop-Artin}, $\mc{X}_\text{pol}$ satisfies Axioms 1,2
and 4 of Corollary~\ref{cor-52}.

The proof of Axiom (3) uses the Grothendieck existence theorem.  Let
$\widehat{A}$ be a complete local algebra and let $(X_n,\mc{L}_n)$ be
a compatible collection of objects of $\mc{X}_\text{pol}$ over $\SP
\widehat{A}/\mf{m}^n$.  For $d$ sufficiently large, $\mc{L}_0^{\otimes
d}$ is very ample and $h^1(X_0,\mc{L}_0^{\otimes d})$ is zero.  The
compatible system of $\widehat{A}/\mf{m}^n$-modules
$H^0(X_n,\mc{L}_n^{\otimes d})$ defines a finite free
$\widehat{A}$-module.  Choosing a basis for this module, there are
induced closed immersions $\mc{X}_n \rightarrow
\PP^N_{\widehat{A}/\mf{m}^n}$ such that the pullback of
$\OO_{\PP^N}(1)$ is $\mc{L}_n^{\otimes d}$.  By the Grothendieck
existence theorem, ~\cite[Corollaire 5.1.8]{EGA3}, there exists a
closed subscheme $X$ of $\PP^N_{\widehat{A}}$ whose reductions give
the compatible family $\mc{X}_n$.  By the theorem on formal functions,
~\cite[Th\'eor\`eme 4.1.5]{EGA3}, $H^0(X,\mc{L}^{\otimes e})$ equals
$\varprojlim H^0(X_n,\mc{L}_n^{\otimes e})$ for every $e$.  Because
the higher cohomologies vanish, and because the $X_n$ are flat over
$\SP \widehat{A}/\mf{m}^n$, also $H^0(X,\mc{L}^{\otimes e})$ is flat
over $\widehat{A}$.  Therefore $X$ is flat over $\SP \widehat{A}$.
Finally, by ~\cite[Corollaire 5.1.6]{EGA3}, there exists an invertible
sheaf $\mc{L}$ on $X$ whose reductions give $\mc{L}_n$.  By
~\cite[Th\'eor\`eme 4.7.1]{EGA3}, this is ample.
\end{proof}

\begin{rmk} \label{rmk-pol}
There are, of course, other ways to verify the proposition.  One can
use the existence of Hilbert schemes to give a smooth morphism from a
scheme to $\mc{X}_\text{pol}$.
\end{rmk}

Let $S$ be an excellent scheme.  Let $\mc{Y}$ be an limit preserving,
separated algebraic stack over $(\text{Aff}/S)$ whose diagonal is
representable by finite morphisms.  Define $\mc{H}$ to be the category
whose objects are 4-tuples $(U,X,L,g)$ where $(U,X)$ is an object of
$\mc{X}_\pflpf$, $L$ is an invertible sheaf on $X$, and
$g:X\rightarrow \mc{Y}$ is a 1-morphism.  A datum is required to
satisfy the condition that $L$ is $g$-relatively ample, i.e., for
every affine scheme $\SP A$ and object $a$ of $\mc{Y}(\text{SP} A)$,
the pullback of $L$ to the 2-fibered product $\SP A\times_{a,\mc{Y},g}
X$ is ample.

Morphisms in $\mc{H}$ are data $(f_U,f_X,f_L,f_g)$ of a morphism
$(f_U,f_X,f_L)$ in $\text{Coh}_\pi$ together with a 2-isomorphism from
the 1-morphism $g':X'\rightarrow \mc{Y}$ to the composite 1-morphism
$$
X'\xrightarrow{f_X} X''\times_{U''} U' \xrightarrow{\text{pr}_{X''}}
X'' \xrightarrow{g''} \mc{Y}.
$$

\begin{prop} \label{prop-Vistoli}
The category $\mc{H}$ is a limit preserving algebraic stack over
$(\text{Aff}/S)$ with quasi-compact, separated diagonal.
\end{prop}

\begin{proof}
There is a 1-morphism $\text{T}:\mc{H}\rightarrow \mc{Coh}_\pi$.  In
fact this factors through the open substack of data $(U,X,L)$ such
that $L$ is an invertible sheaf.  By ~\cite[Theorem 1.1]{OHom}, and
using ~\cite[Th\'eor\`eme 4.7.1]{EGA3} (locally over the Hom stack,
for a quasi-compact smooth cover of $\mc{Y}$ whose image contains the
local image of $g$ ), the 1-morphism $T$ is representable by limit
preserving algebraic stacks.  Moreover, the diagonal morphism
associated to $T$ is quasi-compact and separated.  By Proposition
~\ref{prop-Pic}, the composite $G\circ F:\mc{H} \rightarrow
\mc{X}_\pflpf$ is representable by limit preserving algebraic stacks
with quasi-compact, separated diagonal.

By Proposition ~\ref{prop-comp1}, Proposition ~\ref{prop-comp2},
Corollary ~\ref{cor-algobs}, Proposition~\ref{prop-pf} and
Proposition~\ref{prop-Artin}, $\mc{H}$ satisfies Axioms 1,2 and 4 of
Corollary ~\ref{cor-52}.  It only remains to verify Axiom (3).

To prove Axiom (3), apply the straightforward analogue of the argument
from Proposition~\ref{prop-pol}.  Namely, given a compatible family
$(X_n,L_n,g_n)$ over $\SP \widehat{A}/\mf{m}^n$, first form the
sequence of sheaves $E_{d,n} = (g_n)_*L_n^{\otimes d}$ on $\SP
\widehat{A}/\mf{m}^n \times_S \mc{Y}$.  For $d$ sufficiently large,
the natural map $g_n^*E_{n,d} \rightarrow L_n^{\otimes d}$ is
surjective for every $d$.  For $d$ perhaps larger, the induced
morphism $X_n \rightarrow \text{Proj}\text{Sym}(E_{n,d})$ is a closed
immersion.  These statements are local on $\mc{Y}$, and thus can be
checked after base change by a smooth morphism $\SP B \rightarrow
\mc{Y}$ whose image contains $g_n(X_n)$.  Then the statements follow
from the usual versions for schemes.

The sheaves $E_{n,d}$ are coherent with proper support.  Therefore, by
the analogue of the Grothendieck existence theorem for stacks,
~\cite[Proposition 2.1]{OS}, there exists a coherent sheaf with proper
support $E_d$ on $\SP\widehat{A} \times_S \mc{Y}$ whose reductions are
the sheaves $E_{n,d}$.  Denote by $P$ the stack
$\text{Proj}\text{Sym}(E_d)$ (which can be constructed by flat
descent, for example).  The projection,
$$
\text{pr}: P \rightarrow \SP
\widehat{A}\times_S \mc{Y}, 
$$
is representable by proper, finitely presented algebraic spaces.
Therefore $P$ is a limit preserving algebraic stack with
quasi-compact, separated diagonals.  There are natural closed
immersions $X_n \rightarrow \SP \widehat{A}/\mf{m}^n \times_{\SP
\widehat{A}} P$.  Again by the Grothendieck existence theorem
~\cite[Proposition 2.1]{OS}, there exists a closed substack $X$ of $P$
whose reductions are $X_n$.  Since the reductions are proper, flat,
finitely presented algebraic spaces over $\SP \widehat{A}/\mf{m}^n$,
by the same argument as in the proof of Proposition ~\ref{prop-pol},
$X$ is a proper, flat, finitely presented algebraic space over $\SP
\widehat{A}$.  Define $g$ to be the restriction of $\text{pr}$ to $X$.
By the Grothendieck existence theorem for algebraic spaces,
~\cite[Theorem V.6.3]{Kn}, there exists an invertible sheaf $L$ on $X$
whose reductions are the sheaves $L_n$.  By ~\cite[Th\'eor\`eme
4.7.1]{EGA3} (applied after base change to a quasi-compact smooth
cover of $\mc{Y}$), $L$ is $g$-ample.  Thus $(X,L,g)$ is an object of
$\mc{H}(\SP \widehat{A})$ whose reductions are the objects
$(X_n,L_n,g_n)$.
\end{proof}

\begin{rmk} \label{rmk-final}
There is a stack $\mc{H}_\text{finite}$ closely related to $\mc{H}$
whose objects are data $(U,X,g)$ as above such that $\OO_X$ is
$g$-ample.  There is a 1-morphism $\mc{H}_\text{finite}$ sending
$(U,X,g)$ to $(U,X,\OO_X,g)$.  This morphism is representable by
affine morphisms by the same proof as for the relative
representability of $\text{Coh}_\pi$ in ~\cite[Th\'eor\`eme
4.6.2.1]{LM-B}.  Obviously, if $\OO_X$ is $g$-ample if and only if $g$
is representable by finite morphisms.

Denote by $g:X\rightarrow \mc{H}_\text{finite}\times_S \mc{Y}$ the
universal morphism.  There is an open substack of $X$ where $g$ is
unramified: namely the complement of the support of the sheaf of
relative differentials of $g$.  The complement of this open stack is
proper over $\mc{H}_\text{finite}$.  Its image in
$\mc{H}_\text{finite}$ is proper.  The complement of the image equals
Vistoli's Hilbert stack, cf. ~\cite{VHilb}.

Similarly, using ~\cite[Th\'eor\`eme 12.1.1]{EGA4}, there is an open
substack of $\mc{H}_\text{finite}$ over which the fibers of $X$ are
geometrically reduced and equidimensional.  In the special case that
$\mc{Y}$ is a projective scheme over $S$, this locus is the stack of
\emph{branchvarieties} of Alexeev and Knutson, cf. ~\cite{AlKn}.
\end{rmk}

\bibliography{my}

\begin{thebibliography}{LMB00}

\bibitem[AK06]{AlKn}
Valery Alexeev and Allen Knutson.
\newblock Complete moduli spaces of branchvarieties.
\newblock preprint available at \verb+http://www.math.uga.edu/~valery/+, 2006.

\bibitem[Art69]{Artin}
M.~Artin.
\newblock Algebraization of formal moduli. {I}.
\newblock In {\em Global Analysis (Papers in Honor of K. Kodaira)}, pages
  21--71. Univ. Tokyo Press, Tokyo, 1969.

\bibitem[Art74]{Versal}
M.~Artin.
\newblock Versal deformations and algebraic stacks.
\newblock {\em Invent. Math.}, 27:165--189, 1974.

\bibitem[BGI71]{SGA6}
P.~Berthelot, A.~Grothendieck, and L.~Illusie, editors.
\newblock {\em Th\'eorie des intersections et th\'eor\`eme de
  {R}iemann-{R}och}.
\newblock Springer-Verlag, Berlin, 1971.
\newblock S\'eminaire de G\'eom\'etrie Alg\'ebrique du Bois-Marie 1966--1967
  (SGA 6), Dirig\'e par P. Berthelot, A. Grothendieck et L. Illusie. Avec la
  collaboration de D. Ferrand, J. P. Jouanolou, O. Jussila, S. Kleiman, M.
  Raynaud et J. P. Serre, Lecture Notes in Mathematics, Vol. 225.

\bibitem[CdJ02]{CdJ}
Brian Conrad and A.~J. de~Jong.
\newblock Approximation of versal deformations.
\newblock {\em J. Algebra}, 255(2):489--515, 2002.

\bibitem[Gro63]{EGA3}
A.~Grothendieck.
\newblock \'{E}l\'ements de g\'eom\'etrie alg\'ebrique. {III}. \'{E}tude locale
  des sch\'emas et des morphismes de sch\'emas.
\newblock {\em Inst. Hautes \'Etudes Sci. Publ. Math. 11 (1961), 349-511;
  ibid.}, (17):137--223, 1963.
\newblock \verb+http://www.numdam.org/item?id=PMIHES_1961__11__5_0+.

\bibitem[Gro67]{EGA4}
A.~Grothendieck.
\newblock \'{E}l\'ements de g\'eom\'etrie alg\'ebrique. {IV}. \'{E}tude locale
  des sch\'emas et des morphismes de sch\'emas.
\newblock {\em Inst. Hautes \'Etudes Sci. Publ. Math. 20 (1964), 101-355; ibid.
  24 (1965), 5-231; ibid. 28 (1966), 5-255; ibid.}, (32):5--361, 1967.
\newblock \verb+http://www.numdam.org/item?id=PMIHES_1965__24__5_0+.

\bibitem[Gro03]{SGA1}
A.~Grothendieck.
\newblock {\em Rev\^etements \'etales et groupe fondamental ({SGA} 1)}.
\newblock Documents Math\'ematiques (Paris) [Mathematical Documents (Paris)],
  3. Soci\'et\'e Math\'ematique de France, Paris, 2003.
\newblock S\'eminaire de g\'eom\'etrie alg\'ebrique du Bois Marie 1960--61.
  [Geometric Algebra Seminar of Bois Marie 1960-61], Directed by A.
  Grothendieck, With two papers by M. Raynaud, Updated and annotated reprint of
  the 1971 original [Lecture Notes in Math., 224, Springer, Berlin; MR0354651
  (50 \#7129)].

\bibitem[Ill71]{Illusie1}
Luc Illusie.
\newblock {\em Complexe cotangent et d\'eformations. {I}}.
\newblock Springer-Verlag, Berlin, 1971.
\newblock Lecture Notes in Mathematics, Vol. 239.

\bibitem[Knu71]{Kn}
Donald Knutson.
\newblock {\em Algebraic spaces}.
\newblock Springer-Verlag, Berlin, 1971.
\newblock Lecture Notes in Mathematics, Vol. 203.

\bibitem[LMB00]{LM-B}
G{\'e}rard Laumon and Laurent Moret-Bailly.
\newblock {\em Champs alg\'ebriques}, volume~39 of {\em Ergebnisse der
  Mathematik und ihrer Grenzgebiete. 3. Folge. A Series of Modern Surveys in
  Mathematics [Results in Mathematics and Related Areas. 3rd Series. A Series
  of Modern Surveys in Mathematics]}.
\newblock Springer-Verlag, Berlin, 2000.

\bibitem[Ols04]{OK3}
Martin~C. Olsson.
\newblock Semistable degenerations and period spaces for polarized {$K3$}
  surfaces.
\newblock {\em Duke Math. J.}, 125(1):121--203, 2004.

\bibitem[Ols05a]{ODef}
Martin Olsson.
\newblock Deformation theory of representable morphisms of algebraic stacks.
\newblock preprint, 2005.

\bibitem[Ols05b]{OHom}
Martin Olsson.
\newblock Hom stacks and restriction of scalars.
\newblock preprint, 2005.

\bibitem[OS03]{OS}
Martin Olsson and Jason Starr.
\newblock Quot functors for {D}eligne-{M}umford stacks.
\newblock {\em Comm. Algebra}, 31(8):4069--4096, 2003.
\newblock Special issue in honor of Steven L. Kleiman.

\bibitem[Vis91]{VHilb}
Angelo Vistoli.
\newblock The {H}ilbert stack and the theory of moduli of families.
\newblock In {\em Geometry Seminars, 1988--1991 (Italian) (Bologna,
  1988--1991)}, pages 175--181. Univ. Stud. Bologna, Bologna, 1991.

\end{thebibliography}
\bibliographystyle{alpha}

\end{document}